\documentclass[times, 10pt]{article}
\usepackage{amsmath}
\usepackage{amsthm}

\newtheorem{Theorem}{Theorem}
\newtheorem{Definition}{Definition}
\newtheorem {Lemma} [Theorem]    {Lemma}

\newtheorem {Proposition}[Theorem]    {Proposition}

\newcommand{\by}{{\bf y}}
\newcommand{\bz}{{\bf z}}
\newcommand{\cm}{{\cal M}}

\begin{document}

\title{Increments of Random Partitions}
\author{\c Serban Nacu\thanks{Department of Statistics, University of California, Berkeley, CA 94720, USA. 
{\tt serban@stat.berkeley.edu}.  
{\tt http://stat-www.berkeley.edu/}$\tilde{\ }${\tt serban}.
Research supported in part by NSF Grant DMS-0071448.}}

\maketitle

\begin{abstract}

For any partition of $\{1, 2, \ldots, n\}$ we define its {\it increments} $X_i, 1 \le i \le n$ by
$X_i = 1$ if $i$ is the smallest element in the partition block that contains it, $X_i = 0$
otherwise. We prove that for partially exchangeable random partitions
(where the probability of a partition depends only on its block sizes in order of appearance),
the law of the increments uniquely determines the law of the partition. 
One consequence is that the Chinese Restaurant Process CRP($\theta$)
(the partition with distribution given by the Ewens sampling formula with parameter $\theta$)
is the only exchangeable random partition with independent increments.
\end{abstract}

\section{Introduction}

Random partitions have been studied extensively during the past thirty years, and have found
various applications in population biology, Bayesian statistics, combinatorics, and statistical
physics; see~\cite{JP2} for an in-depth survey. Exchangeable random partitions 
were introduced by Kingman, motivated by applications in genetics. 
Partially exchangeable random partitions were introduced by Pitman in~\cite{JP1}.
We recall their definition:

\begin{Definition}
Consider a random partition $\Pi_n = \{A_1, A_2, \ldots, A_k\}$ of $[n] = \{1, 2, \ldots, n\}$, where the blocks $A_i$ 
are listed in order of appearance (i.e. in increasing order of their smallest element). 
$\Pi_n$ is called {\bf partially exchangeable} if its probability only depends on the ordered block sizes:

\begin{equation}\label{pep}
P(\Pi_n) = p(|A_1|, \ldots, |A_k|)
\end{equation}

\noindent for some function $p$ taking values on the set of {\it compositions} of $n$ (ordered sets of positive integers that sum up to $n$). 
$\Pi_n$ is called {\bf exchangeable} if $p$ is symmetric, so its probability depends on the sizes of the blocks but not their order.
\end{Definition}

A natural way to look at a partition is to construct it one element at a time. We start with a single block $\{1\}$, then the next element $2$ either joins the existing block or starts a new one, and so on until $n$. For $i=1,\ldots, n$ let $X_i = 1$ if $i$ starts a new block, $X_i = 0$ otherwise. (Alternatively, $X_i=1$ if $i$ is the smallest element in its block, $X_i=0$ otherwise.) We always define $X_1=1$. 

\begin{Definition} Let $\Pi_n$ be a partition of $\{1, 2, \ldots, n\}$.
Let $X_i = 1$ if $i$ is the smallest element in its block, $0$ otherwise.
We call the sequence $(X_1, \ldots, X_n)$ the {\bf increments} of $\Pi_n$.
\end{Definition}

Clearly the partition determines the increments, but not viceversa. The increments do not even determine the block sizes: the partitions $\{\{1,3,4\},\{2\}\}$ and $\{\{1,3\},\{2,4\}\}$ have the same increments $X_1=1, X_2=1, X_3=0, X_4=0$.

For a random partition, the law of the partition induces a law for its increments on the probability space $\{0,1\}^n$. We prove that if the partition is partially exchangeable, the law of the increments {\it does} determine the law of the partition:

\begin{Theorem}\label{main}
If $\Pi_n$ is a partially exchangeable partition of $[n]$, the distribution of its increments $(X_1, \ldots, X_n)$ uniquely determines the distribution of $\Pi_n$.
\end{Theorem}

Hence we obtain a correspondence between the set of partially exchangeable laws for partitions, and the set of laws for binary sequences of zeroes and ones. This correspondence is as close to a bijection as we could possibly hope for, in the following sense. There are exactly $2^{n-1}$ compositions of $n$, so according to~(\ref{pep}), the law of the partition is determined by the $2^{n-1}$ values of $p$, which are arbitrary except for the constraint that $P$ is a probability measure (so $p \ge 0$ and a weighted sum of its values is equal to 1). There are $2^{n-1}$ possible sequences of increments (since $X_1=1$ always), so their law is also determined by $2^{n-1}$ numbers; the constraint on them turns out to be a system of linear inequalities (plus the obvious constraint that they sum up to 1).

As an application of Theorem~\ref{main}, we answer a question raised by Jim Pitman. First we define
\begin{Definition}
Let $\theta > 0$. The {\bf Chinese Restaurant Process} with parameter $\theta$ is the exchangeable
random partition of $[n]$ with distribution given by
\begin{equation}\label{crpdist}
p(n_1, \ldots, n_k) = \theta^k \prod_{i=1}^k (n_i-1)! / \prod_{i=0}^{n-1} (\theta + i)
\end{equation}
We denote it by $CRP(\theta)$.
\end{Definition}

Equation~(\ref{crpdist}) is equivalent to the well-known Ewens sampling formula, and the distribution
of the block sizes of $CRP(\theta)$ is also referred to as the Ewens partition structure. It has also
been referred to as the Blackwell-MacQueen distribution. See~\cite{TE} for a survey. The name
"Chinese Restaurant Process" was introduced by Lester Dubins and Jim Pitman in the early 80's.

$CRP(\theta)$ has several equivalent descriptions; for example, for $\theta = 1$ it is the partition
induced by the cycles of an uniform random permutation. See~\cite{JP2} for details.
The description we are interested in is in
terms of its increments $X_i$. It follows from~(\ref{crpdist}) that those are independent and satisfy
\begin{equation}
P(X_i = 1) = \theta / (i - 1 + \theta)
\end{equation}

Hence $CRP(\theta)$ admits a simple construction one element at a time: $i$ starts a new block with
probability $\theta / (i - 1 + \theta)$ (and joins an existing block with probability proportional to the size
of the block). It is easy to prove that this constructs an exchangeable partition.
Jim Pitman asked whether this is the only exchangeable random partition with the
property that $X_i$ are independent. We prove that the answer is yes:

\begin{Theorem}\label{CRP}
The Chinese Restaurant Process $CRP(\theta)$ is the only exchangeable random partition with independent increments.
\end{Theorem}

Hence in this sense, the Chinese Restaurant Process is the simplest exchangeable random partition. \\

Several other random partitions admit simple representations in terms of their increments $X_i$. If the
increments are not independent, the next simplest case is to assume some kind of Markov structure.
For example, we can require that the partial sums $S_i = X_1 + \ldots + X_i$ form a Markov chain.
This is the same as requiring that the probability that $i+1$ start a new block depend only on $i$ and
on the number of already existing blocks $S_i$. One process that satisfies this is
Pitman's two-parameter generalization of CRP($\theta$), described in~\cite{JP1},~\cite{JP2}.
In this case we have
$$P(X_{i+1}=1 | S_i = k) = (k \alpha + \theta) / (i + \theta)$$
where the parameters satisfy $0 \le \alpha < 1, \theta > -\alpha$. For $\alpha = 0$ we obtain
CRP($\theta$). It is an open question to describe all random partitions for which $S_i$ is a Markov chain.

\section{Proof of Theorem~\ref{main}}

Let $\Pi_n = \{A_1, A_2, \ldots, A_k\}$ be a partition of $[n]$ with $k$ blocks. We list the blocks $A_i$ in order of appearance, so $1 \in A_1$, the smallest element not in $A_1$ is in $A_2$, and so on. Let $B(\Pi_n) = (|A_1|, \ldots, |A_k|)$ the sizes of its blocks in order of appearance.

Let $(X_1, \ldots, X_n)$ be the increments of the partition. Since there are $k$ blocks, there are exactly $k$ ones and $n-k$ zeroes among the increments. We can encode such a binary sequence by the distance between consecutive 1's. For example, the sequence $1,0,0,1,1,0,1,0$ will be encoded $3, 1, 2, 2$ as the distance between the first and the second 1 is three, the distance between the second and the third 1 is one, and so on. Formally, if $a_i$ is the smallest element in $A_i$, then we know $X_{a_i} = 1$ for all $1 \le i \le k$, so we encode the increments as the $k$-tuple $D(\Pi_n) = (a_2 - a_1, a_3 - a_2, \ldots, a_k - a_{k-1}, n+1-a_k)$. There is a bijection between such $k$-tuples and binary sequences with $X_1 = 1$, so from now on we will identify the sequence of increments with its encoding, and work directly with $k$-tuples.

Let $S_{n,k}$ be the set of $k$-tuples whose elements add up to $n$. We are interested in the relationship between $B(\Pi_n)$ and $D(\Pi_n)$; both are in $S_{n,k}$. We define a partial order relation on $S_{n,k}$ as follows: $(y_1, \ldots, y_n) \ge (z_1, \ldots, z_n)$ iff $y_1 \ge z_1, y_1 + y_2 \ge z_1 + z_2, \ldots, y_1+y_2+\ldots+y_n \ge z_1+z_2+\ldots+z_n$. Then we have

\begin{Lemma}\label{BD}
For any partition $\Pi_n$, $B(\Pi_n) \ge D(\Pi_n)$. 
\end{Lemma}

\begin{proof} Let $B(\Pi_n) = (b_1, \ldots, b_k), D(\Pi_n) = (d_1, \ldots, d_k)$. Fix $m$ with $1 \le m \le k$. Clearly $\sum_{i=1}^m b_i = |\bigcup_{i=1}^m A_i|$ and $\sum_{i=1}^m d_i = a_{m+1} - a_1 = a_{m+1} - 1$, where $a_{m+1}$ is the smallest element in $A_{m+1}$. Since the blocks $A_i$ are listed in order of appearance, $a_{m+1}$ is the smallest element outside $\bigcup_{i=1}^m A_i$, so it is at most $ |\bigcup_{i=1}^m A_i| + 1$; equality occurs iff  $\bigcup_{i=1}^m A_i = \{1, 2, \ldots, a_{m+1}-1\}$.
\end{proof}

Now consider a partially exchangeable law for $\Pi_n$, defined as in~(\ref{pep}) by a function $p$. This induces a law for the increments of $\Pi_n$, and by using the encoding of binary sequences into $k$-tuples discussed above, this induces a law on $k$-tuples:

\begin{equation}\label{qeq}
q(d_1, \ldots, d_k) = P(D(\Pi_n) = (d_1, \ldots, d_k))
\end{equation}

\noindent Summing up over all partitions with the same block structure, we obtain

\begin{equation}
q(d_1, \ldots, d_k) = \sum_{(b_1, \ldots, b_k) \in S_{n,k}} p(b_1, \ldots, b_k) r(d_1, \ldots, d_k; b_1, \ldots, b_k)
\end{equation}

\noindent where $r(a_1, \ldots, a_k; b_1, \ldots, b_k)$ denotes the number of partitions $\Pi_n$ with blocks $B(\Pi_n) = (b_1, \ldots, b_k)$ and increments encoded as $D(\Pi_n) = (d_1, \ldots, d_k)$.

This gives a system of linear equations in the $q$'s and the $p$'s; we will show it can be solved to compute $p$ in terms of $q$. Consider the dictionary order on $S_{n,k}$: 
$(y_1, \ldots, y_k) \ge_d (z_1, \ldots, z_k)$ iff $y_1 > z_1$ or $y_1 = z_1, \ldots, y_m = z_m$ and $y_{m+1} > z_{m+1}$ for some $m$ or $y_i = z_i$ for all $i$. It is easy to see that if $(y_1, \ldots, y_k) \ge (z_1, \ldots, z_k)$ in the partial order previously defined, then $(y_1, \ldots, y_k) \ge_d (z_1, \ldots, z_k)$ in the dictionary order. Hence from Lemma~\ref{BD} we obtain

\begin{Lemma}
For $\by, \bz \in S_{n,k}$, $r(\by; \bz) = 0$ unless $\bz \ge \by$, and $r(\by; \bz) = 1$ if $\bz = \by$.
\end{Lemma}

Now $S_{n,k}$ is totally ordered under the dictionary order so we can arrange its elements in decreasing order. But then the lemma says that the matrix of the system of linear equations is triangular and its diagonal elements are all 1, so it is trivially invertible and $p$ can be computed in terms of $q$. Explicitly, if the elements of $S_{n,k}$ are ${\bf y_1} >_d {\bf y_2} >_d {\bf y_3} > \ldots$, then $p({\bf y_1}) = q({\bf y_1})$ and for $i \ge 2, p({\bf y_i}) = q({\bf y_i}) - \sum_{j=1}^{i-1} p({\bf y_j}) r({\bf y_i};{\bf y_j})$. This completes the proof of Theorem~\ref{main}. \qed \\

The following result gives an explicit formula and a generating function for the coefficients $r(\cdot; \cdot)$. 

\begin{Proposition}
For $(d_1, \ldots, d_k) \in S_{n,k}$ and $(b_1, \ldots, b_k) \in S_{n,k}$,
let $\cm$ be the set of $k \times k$ square matrices $M = (m_{ij})$ with the following properties: 
\begin{itemize}
\item[(i)] $M$ is upper triangular, so $m_{ij} = 0$ if $i>j$.
\item[(ii)] All other entries are non-negative integers, so $m_{ij} \ge 0$ if $i \le j$.
\item[(iii)] For all $i$, the $i$-th row sums up to $b_i - 1$.
\item[(iv)] For all $i$, the $i$-th column sums up to $d_i - 1$.
\end{itemize}

Then
\begin{equation}
r(d_1, \ldots, d_k; b_1, \ldots, b_k) = \sum_{M \in \cm} \prod_i (d_i - 1)! / \prod_{i,j} (m_{ij})!
\end{equation}

Hence the following generating function identity holds:
\begin{eqnarray}
\lefteqn{x_1^{d_1-1} (x_1+x_2)^{d_2-1} \ldots (x_1+x_2+\ldots+x_k)^{d_k-1} =} \quad \nonumber \\
& & \sum_{{\bf b} \in S_{n,k}} r(d_1, \ldots, d_k; b_1, \ldots, b_k) x_1^{b_1-1} x_2^{b_2-1} \ldots x_k^{b_k-1}
\end{eqnarray}

\end{Proposition}

\begin{proof}
We need to count the number of partitions with block structure ${\bf b}$ and increment structure ${\bf d}$. Let $A_i$ be the blocks of such a partition, in order of appearance. Knowing ${\bf d}$ is equivalent to knowing the smallest element of each $A_i$, call them $a_i$. Let $D_i = \{a_i, a_i +1, \ldots, a_{i+1} - 1\}$ and let $n_{ij} = |A_i \bigcap D_j|$. Of the elements in $D_j$, $a_j$ must belong to $A_j$; all others may be assigned to any $A_i$ with $i \le j$, and the number of ways in which this can be done is $(d_j - 1)! / (n_{1j}! n_{2j}! \ldots n_{j-1 \, j}! (n_{jj}-1)! )$. If we let $m_{ij} = n_{ij}$ if $i \ne j$, $m_{ii} = n_{ii}-1$, then we obtain the desired formula. The generating function identity follows easily.
\end{proof}

\section{Partitions with Independent Increments}

We now prove Theorem~\ref{CRP}.

\begin{proof} Let $\Pi_n$ be an exchangeable random partition of $[n]$ such that its increments $X_1, \ldots, X_n$ are independent. 
As above, let $p$ be the function that describes its distribution as in~(\ref{pep}), and $q$ the joint distribution of the increments.
Consider sequences of increments with $k=n-1$; that is, there is only one zero among $X_1, \ldots, X_n$. 
If the zero is at the beginning ($X_2=0$) then the partition must be $\Pi_n=\{\{1,2\}, \{3\},\ldots, \{n\}\}$. Hence
\begin{equation*}
q(1,0,1,\ldots,1)=p(2,1,\ldots,1). 
\end{equation*}

If the zero is at the end ($X_n = 0$) then $n$ could be in any of the $n-1$ pre-existing blocks so there are $n-1$ choices for $\Pi_n$. Since $\Pi_n$ is exchangeable, all these choices are equally likely. Hence 
$$q(1,1,1,\ldots,1, 0)=(n-1)p(2,1,\ldots,1).$$

But if we let $u_n = P(X_n=1)$, by independence we also have 
$$q(1,0,1,\ldots,1)=u_1 (1-u_2) u_3 \ldots u_n$$ 
and 
$$q(1,1,1,\ldots,1, 0)=u_1 u_2 \ldots u_{n-1} (1-u_n).$$ 
If $u_2 \ne 0,1$ then it follows easily by induction that $u_n \ne 0$ for all $n$ and hence 
$$(n-1)(1-u_2)/u_2 = (1-u_n)/u_n$$
so if we let $\theta=u_2/(1-u_2)$ then $u_n = \theta / (n-1+\theta)$. Hence the increments of the process have the same law as the increments of CRP($\theta$), so the process must be CRP($\theta$).

It remains to consider the cases when $u_2$ is 0 or 1. If $u_2=1$ then by induction $u_n=1$ so all blocks are singletons (this is CRP($\theta$) in the limit case $\theta=\infty$). If $u_2=0$ then all $u_n=0$ ($\Pi_n$ cannot contain the singleton \{1\} so by exchangeability it cannot contain the singleton $\{n\}$ either) so there is only one block (CRP($\theta$) for $\theta=0$).
\end{proof}

\section{Another Binary Representation}

Theorem~\ref{main} allows us to obtain a partially exchangeable partition from any law on random binary sequences satisfying certain constraints. While we have proved the theorem for partitions of the finite set $[n]$, the result
is easily extended to infinite partitions. We obtain thus a correspondence between the set of distributions of partially exchangeable partitions of ${\bf N}$, and the set of distributions of infinite binary sequences. By Theorem~\ref{main}, the correspondence is one-to-one; it is not onto.

There is another way to associate infinite binary sequences to partitions, which is discussed in detail in~\cite{Y}, Chapter 4. Given a binary sequence, the problem is to construct an exchangeable random partition $\Pi$, so that the distribution of the (unordered) block sizes of its restriction $\Pi_n$ to $[n]$ is the same as the distribution of the ``gaps'' (distances) between consecutive 1's in the binary sequence. More precisely:

\begin{Definition}
Let $Y_1, Y_2, \ldots$ be a random infinite sequence with $Y_1 = 1$ and $Y_n \in \{0,1\} \, \forall n$, and let $\Pi$ be
an exchangeable random partition of ${\bf N}$. Let $1 = n_1 < n_2 < \ldots$ be the locations of the $1$'s in the
sequence $Y_n$. For any fixed $n$, let $n_k = \max \{i: i \le n, X_i = 1\}$. Then
$$ n = (n_2 - n_1) + \ldots + (n_k - n_{k+1}) + (n+1-n_k) $$
is a partition of $n$ into $k$ integers. If this has the same distribution
as the partition of $n$ induced by restricting $\Pi$ to $[n]$, then we say that
$\Pi$ has a {\bf gap representation} by $Y$.
\end{Definition}

The two binary representations are related, but not identical. In particular, while any random partition has a representation by increments, it is not known under what conditions a random partition admits a gap representation. Also, note that the gap representation is interesting only for infinite partitions and sequences; for fixed $n$, any law for the sequence easily translates into a law for an exchangeable partition, as we are free to specify to probabilities for all possible block sizes. The problem is whether these laws are compatible as $n$ varies.

In~\cite{Y}, gap representations are constructed for various partitions, including CRP($\theta$), for
which the gap representation has $P(Y_n = 1) = \theta / (n - 1 + \theta)$ and the $Y_n$ are
independent. It is also proven that
the only exchangeable random partition which admits a gap representation via a sequence $Y$ of independent
binary random variables is CRP($\theta$). The similarity with Theorem~\ref{CRP} may seem surprising, but it is
explained by the following result:

\begin{Proposition}
Suppose $\Pi$ has a gap representation by $Y$, and let $X$ be the increments of $\Pi$. Then
\begin{equation}\label{xeqy}
X_1 + \ldots + X_n \stackrel{d}{=} Y_1 + \ldots + Y_n, \quad \forall n \ge 1
\end{equation}
\end{Proposition}

\begin{proof}
Both sides of the identity are equal in distribution to the number of blocks in $\Pi_n$.
\end{proof}

If $X$ and $Y$ are each sequences of independent variables, then (\ref{xeqy}) implies
they have the same distribution (in fact, it is enough to assume that the sequences of
partial sums $X_1 + \ldots + X_n$ and $Y_1 + \ldots + Y_n$ are 
(possibly non-homogenous) Markov chains). \\

{\bf Acknowledgments.} Many thanks to Jim Pitman for suggesting the problem
and for useful discussions. Thanks also to Noam Berger.

\end{document}